\documentclass[reqno]{amsart}
\usepackage{amssymb}
\usepackage{hyperref}

\begin{document}
\title[ ]
{ On  $\lambda$-commuting operators}
\author[ A. Tajmouati,  A. El Bakkali and M.B. Mohamed Ahmed ]
{  A. Tajmouati, A. El Bakkali and M.B. Mohamed Ahmed}

\address{A. Tajmouati and M.B. Mohamed Ahmed \newline
 Sidi Mohamed Ben Abdellah
 Univeristy
 Faculty of Sciences Dhar Al Mahraz Fez, Morocco.}
\email{abdelaziztajmouati@yahoo.fr}
\email{bbaba2012@gmail.com}

\address{A. El Bakkali  \newline
Department of Mathematics
University Chouaib Doukkali,
Faculty of Sciences, Eljadida.
24000, Eljadida, Morocco.}
\email{aba0101q@yahoo.fr}

\newtheorem{thm}{Theorem}[section]
\newtheorem{defini}{Definition}[section]
\newtheorem{rem}{Remark}
\newtheorem{lem}{Lemma}[section]
\newtheorem{prop}{Proposition}[section]
\newtheorem{coro}{Corollairy}[section]
\newtheorem{ex}{Example}
\newcommand{\pr} {{\bf   Proof: \hspace{0.3cm}}}

\maketitle

\begin{abstract}
In this paper, we study the operator equation $AB=\lambda BA$ for a bounded operator $A,B$ on a complex Hilbert space. We focus on algebraic relations between different operators that include normal, $M$-hyponormal, quasi $*$-paranormal and other classes.

\end{abstract}

\section{Introduction}
Throughout, we will denote by $\mathcal{B}(\mathcal{H})$   the complex Banach algebra of all bounded linear operators  on a   complex  Hilbert space $\mathcal{H}$.
Recall that an operator $T\in \mathcal{B}(\mathcal{H})$ is said to be   :
\begin{itemize}

  \item positive if  $ <Tx,x>\geq 0$ for all $x\in \mathcal{H}$
  \item self-adjoint if  $T=T^*$
  \item isometry if $T^*T=I$, which equivalent to the condition $\|Tx\|=\|x\|$ for all $x\in \mathcal{H}$
  \item normal if $T^*T=TT^*$
  \item unitary $T^*T=TT^*=I$ (i.e $T$ is an onto isometry)
  \item quasinormal if $T[T^*T]=[T^*T]T$
  \item binormal if $[T^*T][TT^*]=[TT^*][T^*T]$ [\ref{3}]
  \item subnormal if $T$ has a normal extension
  \item hyponormal if $T^*T\geq TT^*$, which equivalent to the condition $\|T^*x\|\leq\|Tx\|$ for all $x\in \mathcal{H}$ [\ref{12}]
  \item $M$-hyponormal if $T^*T\geq MTT^*$, where $1\leq M$ which equivalent to the condition $\|T^*x\|\leq M \|Tx\|$ for all $x\in \mathcal{H}$ [\ref{16}]
    \item $p$-hyponormal if $(T^*T)^p\leq (TT^*)^p$, where $0<p\leq 1$ [\ref{1}]
     \item class $A$ if $|T|^2\leq |T^2|$, where $|T|=(T^*T)^{\frac{1}{2}}$
  \item paranormal if $\|Tx\|^2\leq\|T^2x\|\|x\|$ for all $x\in \mathcal{H}$ [\ref{17}]
   \item $k$-hyponormal if $\|Tx\|^k\leq \|T^kx\|\|x\|^{k-1}$ for all $x\in \mathcal{H}$, where $ k\geq 2$
  \item $*$-paranormal if $\|T^*x\|^2\leq\|T^2x\|\|x\|$ for all $x\in \mathcal{H}$ [\ref{18}]
\item quasi $*$-paranormal if $\|T^*Tx\|^2\leq \|T^3x\|\|Tx\|$ for all $x\in \mathcal{H}$   [\ref{19}]
\item $\log$-hyponormal if $T$ invertible and satisfies $\log(T^*T)\geq \log(TT^*)$ [\ref{13}]
  \item $p$-quasihyponormal if $T^*[(T^*T)^p-(TT^*)^p]T\geq 0$ , where $0<p\leq 1$ [\ref{2}]
  \item normoloid if $\|T\| =r(T)$
   \item quasinilpotent if $r(T)= 0$ , where $r(T)=\lim\|T^n\|^{\frac{1}{n}}$
  \end{itemize}

  We can notice that $T$ is hyponormal if $T$ is $p$-hyponormal with $p=1$ and that $p$-hyponormal is $q$-hyponormal for every $0<q\leq p$ from L\"{o}wner-Heinz inquality [\ref{10}].
  Also we can notice that $T$ is paranormal if $T$ is $k$-hyponormal with $k=2$.\\
  It known that invertible $p$-hyponormal is $\log$-hyponormal, we may regard $\log$-hyponormal operator as $0$-hyponormal [\ref{13}].\\

  It is well known that for any operators $A,B$ and $C$, $ A^*A-2\lambda B^*B+\lambda^{2}C^{*}C\geq 0$ for all $\lambda>0$, if and only if $||Bx||^{2}\leq ||Ax||||Cx||\,\mbox{for all}\,x\in H.$
  Thus we have\\
  \begin{itemize}
                 \item $T$ is quasi $*$-paranormal if and only if $T^*[(T^*)^2T^2-2\lambda TT^*+\lambda^2]T\geq 0$, for all $\lambda>0.$
                 \item $T$ is $*$-paranormal if and only if $(T^*)^2T^2-2\lambda TT^*+\lambda^2\geq 0$ , for all $\lambda>0$.\\
               \end{itemize}

  Thus every  $*$-paranormal is quasi $*$-paranormal  and we have the following implications:\\

quasinormal $\Longrightarrow$ binormal \\

  self-adjoint $\Longrightarrow$ normal $\Longrightarrow$ quasinormal $\Longrightarrow$ subnormal $\Longrightarrow$ hyponormal $\Longrightarrow$ $*$-paranormal $\Longrightarrow$ quasi $*$-paranormal\\

  hyponormal $\Longrightarrow$ $p$-hyponormal $\Longrightarrow$ $p$-quasihyponormal $\Longrightarrow$ class $A$ $\Longrightarrow$ paranormal\\

  invertible $p$-hyponormal $\Longrightarrow$ $\log$-hponormal $\Longrightarrow$ paranormal.\\

  For a scalar $\lambda$, two operators $A$ and $B$ in $\mathcal{B}(\mathcal{H})$ are said be $\lambda$-commute if $AB=\lambda BA$. Recently many authors have studied this equation for several classes of operators, for example\\

  \begin{itemize}
    \item In [\ref{9}] authors have proved that if an operator in $\mathcal{B}(\mathcal{H})$ $\lambda$-commutes with a
compact operator, then this operator has a non-trivial hyperinvariant subspace
\item  In [\ref{7}] Conway and Prajitura characterized the closure and the interior of the set of operators that $\lambda$-commute with
a compact operator
\item In [\ref{15}] Zhang ,Ohawada and Cho have studied the properties of operators that  $\lambda$-commute with
a paranormal operator
\item In [\ref{4}] Brooke, Busch and Pearson showed that  if $AB$ is not quasinilpotent, then $|\lambda| =1$,
 and if $A$ or $B$ is self-adjoint then $\lambda\in\mathbb{R}$
 \item In [\ref{14}] Yang and Du gave a simple proofs and generalizations of this results, particulary if $AB$ is bounded below if and only if both $A$ and $B$ are bounded below
 \item In [\ref{11}] Schmoeger  generalized this results
to hermitian or normal elements of a complex Banach algebra
\item In [\ref{5}] Cho, Duggal, Harte and Ota generalized some Schmoeger's results and they
gave the new characterization of a commutativity of Banach space operators.\\
\end{itemize}
The aim of this paper is to study the situation for binormal, $M$-hyponormal, quasi $*$-paranormal operators. Again other related results are also given.\\

We denote the range and the kernel of $T$ by $R(T)$ and $N(T)$ respectively.

\section{ Main results }

We begin with the following result.

\begin{lem}[\ref{19}]\label{lem1}
Let $T \in \mathcal{B}(\mathcal{H})$ be quasi $*$-paranormal and $M$ be a closed $T$-invariant subspace of $\mathcal{H}$. Then  $T_{|M}$ is also  quasi $*$-paranormal.
\end{lem}

\begin{lem}
Let $A\in \mathcal{B}(\mathcal{H})$ be quasi $*$-paranormal operator. If  $A$ is quasinilpotent, then $A=0$
\end{lem}

\begin{pr}
Suppose that $A\neq0$, we have for any operator $  T\in \mathcal{B}(\mathcal{H})$,
$\|T^*T\|=\|TT^*\|=\|T\|^2=\|T^*\|^2$.
Firstly, by induction we prove that $A$ is normaloid.\\

For $n=3$,
since $A\in \mathcal{B}(\mathcal{H})$ is quasi $*$-paranormal, then
$\|A\|^4=[\|A\|^2]^2=\|A^*A\|^2\leq\|A^3\|\|A\|$, whence

$\|A\|^3\leq\|A^3\|\leq \|A\|^3$. From where $\|A^3\|= \|A\|^3$.\\
It is now assumed  for an integer $n\geq 3$ and any integer k such that $3\leq k\leq n$,
we have $\|A^k\|= \|A\|^k$.
we  know that $\|A^{n+1}\|\leq\|A\|^{n+1}$, other hand We have
\begin{eqnarray*}
[\|A\|^{n-1}]^4  & =& [\|A^{n-1}\|^2]^2\\
                 & = & [\sup_{\|x\|=1}<A^{n-1}x,A^{n-1}x>]^2\\
                  & = & [\sup_{\|x\|=1}<A^*A[A^{n-2}x],A^{n-2}x>]^2 \\
                  & \leq & \sup_{\|x\|=1} [\|A^*A[A^{n-2}x]\|]^2[\|A^{n-2}x\|]^2\\
                   & \leq & \sup_{\|x\|=1} [\|A^3[A^{n-2}x]\|\|A[A^{n-2}x]\|[\|A^{n-2}x\|]^2]\\
                   & \leq& \sup_{\|x\|=1} [\|A^{n+1}x]\|\|A^{n-1}x\|[\|A^{n-2}x\|]^2]\\
                    & \leq & \sup_{\|x\|=1} [\|A^{n+1}x]\|\sup_{\|x\|=1}\|A^{n-1}x\|\sup_{\|x\|=1}\|A^{n-2}x\|^2]\\
                     & \leq & \|A^{n+1}\|\|A^{n-1}\|\|A^{n-2}\|^2\\
                        & \leq & \|A^{n+1}\|\|A\|^{n-1}\|A\|^{2(n-2)}
\end{eqnarray*}
whence $\|A^{n+1}\| \geq  \|A\|^{4(n-1)-(n-1)-2(n-2)}=\|A\|^{n+1}$,
thus   $\|A^{n+1}\|=\|A\|^{n+1}$. Since $A$ is normaloid we have $r(A)=\|A\|$.\\
But $A$ is quasinilpotent, then
$r(A)=0$, therefore  $\|A\|=0$  and this is contradiction.
\end{pr}

\begin{coro}
Let $A\in \mathcal{B}(\mathcal{H})$ be $*$-paranormal operator. If  $A$ is quasinilpotent, then $A=0$
\end{coro}

\begin{thm}
 Let $A,B\in \mathcal{B}(\mathcal{H})$ satisfy $AB=\lambda BA\neq 0,\> \lambda\in\mathbb{C}$. Let $A$ be hyponormal and $B$ be quasi $*$-paranormal. If $A$ is invertible or 0 is an isolated point of $\sigma(A)$, then $|\lambda|=1$.
\end{thm}

\begin{pr}
It is clear that if $A$ is invertible then $|\lambda|=1$.\\
Suppose that 0 is  isolated point of $\sigma(A)$.
Since  $\mathcal{H}=\overline{R(A^*)} \oplus N(A)$  we can decompose $A=\left(
               \begin{array}{cc}
                 A_1 & 0 \\
                 A_2 & 0 \\
               \end{array}
             \right)$
on $\overline{R(A^*)} \oplus N(A)$.\\
Since $A$ is hyponormal
$$\left(\begin{array}{cc}
          A_1^*A_ 1+A_2^*A_2& 0 \\
          0 & 0 \\
        \end{array}
      \right)=A^*A\geq
      AA^*=\left(
              \begin{array}{cc}
                A_1A_1^* & A_1A_2^* \\
                A_2A_1^* & A_2A_2^* \\
              \end{array}
            \right)$$
it holds $ A_2A_2^*\leq 0$
and hence $A_2 =0$, we have $$A=\left(
                        \begin{array}{cc}
                          A_1 & 0 \\
                          0 & 0 \\
                        \end{array}
                      \right),
$$
Consequently  $\sigma(A)=\sigma(A_1)\cup\{0\}$.\\
We conclude that if $A_1$ is not invertible then 0 is isolated point of $\sigma(A_1)$ and since $A_1$ is hyponormal, then  $0\in\sigma_p(A_1)$, a contradiction and hence $A_1$ is invertible.\\
Let $B=\left(
                                                     \begin{array}{cc}
                                                       B_1 & B_2 \\
                                                       B_3 & B_4 \\
                                                     \end{array}
                                                   \right)$
                                                   on $\overline{R(A^*)} \oplus N(A)$.
     By $\lambda\neq 0$ and $AB=\lambda BA$, it holds $A_1B_2=B_3A_1=0$. Since $A_1$ is invertible, we have $B_2=B_3=0$. Hence, $$B=\left(
                                               \begin{array}{cc}
                                                 B_1 & 0 \\
                                                 0 & B_4 \\
                                               \end{array}
                                             \right).$$
                                              Therefore $\overline{R(A^*)}$ is invariant for $B$ and hence $B_1=B_{|\overline{R(A^*)}}$ is quasi $*$-paranormal by lemma$\>$\ref{lem1}.
                                               Since
                                               \begin{eqnarray*}
                                                AB &=&   \left(\begin{array}{cc}  A_1 & 0 \\ 0 & 0 \\
                        \end{array}
                      \right)\left(\begin{array}{cc}
                                                 B_1 & 0 \\
                                                 0 & B_4 \\
                                               \end{array}
                                             \right)\\
                                             &=&
                                         \left(\begin{array}{cc}
                                                 A_1B_1 & 0 \\
                                                 0 & 0\\
                                               \end{array}
                                             \right)\\
                                              &=&
                                             \lambda BA\\
                                             &= & \lambda \left(\begin{array}{cc}
                                                 B_1 & 0 \\
                                                 0 & B_4 \\
                                               \end{array}
                                             \right)\left(
                        \begin{array}{cc}
                          A_1 & 0 \\
                          0 & 0 \\
                        \end{array}
                      \right)\\
                        &= & \lambda  \left(\begin{array}{cc}
                                                 B_1A_1 & 0 \\
                                                 0 & 0 \\
                                               \end{array}
                                             \right).
                                              \end{eqnarray*}

                                             Whence  $A_1B_1=\lambda B_1A_1$, then $B_1=\lambda A_1^{-1}B_1A_1$.\\
                                             If $r(B_1)= 0$ then $B_1$ is quasi $*$-paranormal and quasinilpotent, then $B_1=0$ by lemma 2.2, therefore AB=0 and it's a  contradiction  because $AB\neq 0$. From where $r(B_1)\neq 0$ and consequently  $|\lambda|=1$.
\end{pr}

\begin{coro}
 Let $A,B\in \mathcal{B}(\mathcal{H})$ satisfy $AB=\lambda BA\neq 0,\> \lambda\in\mathbb{C}$. Let $A$ be hyponormal and $B$ be $*$-paranormal operator. If $A$ is invertible or 0 is an isolated point of $\sigma(A)$, then $|\lambda|=1$.
\end{coro}

\begin{thm}
Let $A$ be quasinormal operator and $B$ normal such that $AB=\lambda BA\neq 0,\;\; \lambda\in\mathbb{C}$.
If $|\lambda|=1$ then $AB$ is quasinormal operator.
\end{thm}

\begin{pr}
Assume that $AB=\lambda BA\neq 0$, then $B^*A^*=\bar{\lambda}A^*B^*$. Since $B$ and $\lambda B$ are normal operators by Fuglede-Putnam theorem, condition $AB=\lambda BA$ imply that $BA^*=\lambda A^*B$ and $AB^*=\bar{\lambda}B^*A$. Now have
\begin{eqnarray*} AB[(AB)^*AB] &=&[AB][B^*A^*AB]\\
                            &=&[\lambda BA]B^*A^*AB\\
                            &=&\lambda B[AB^*]A^*AB\\
                                 &=&\lambda B[\bar{\lambda}B^*A]A^*AB\\
                                 &=&|\lambda|^2 [BB^*][AA^*A]B\\
                                 &=& [B^*B][A^*AA]B\\
                                  &=& B^*[BA^*]AAB\\
                                     &=& B^*[\lambda A^*B]AAB\\
                                     &=& B^*A^*[\lambda BA]AB\\
                                     &=& B^*A^*[AB]AB\\
                                      &=& [(AB)^*AB] AB
                                      \end{eqnarray*}
                                      Therefore  $AB$ is quasinormal.\\
\end{pr}
\begin{thm}
Let $A$ be binormal operator and $B$ normal such that $AB=\lambda BA\neq 0, \;\; \lambda\in\mathbb{C}$
if $|\lambda|=1$ then $AB$ is binormal operator.
\end{thm}

\begin{pr}
Because $B$ and $\lambda B$ are normal operators by Fuglede-Putnam theorem, condition $AB=\lambda BA$ imply that
$BA^*=\lambda A^*B$ and $AB^*=\bar{\lambda}B^*A$. Therfore
\begin{eqnarray*}AB(AB)^*(AB)^*AB &=& A[BB^*]A^*B^*A^*AB\\
                 &=& A[B^*B]A^*B^*A^*AB\\
                  &=& [AB^*]BA^*[B^*A^*]AB\\
                   &=& [\bar{\lambda}B^*A]BA^*[\bar{\lambda}A^*B^*]AB\\
                     &=& (\bar{\lambda})^2B^*[AB]A^*A^*[B^*A]B\\
                          &=& (\bar{\lambda})^2B^*[\lambda BA]A^*A^*[\frac{1}{\bar{\lambda}}AB^*]B\\
                   &=& |\lambda|^2B^*B[AA^*A^*A]B^*B\\
                    &=& B^*B[A^*AAA^*]B^*B\\
                     &=& B^*[BA^*]AA[A^*B^*]B\\
                     &=& B^*[\lambda A^*B]AA[\frac{1}{\bar{\lambda}}B^*A^*]B\\
                     &=& \frac{\lambda}{\bar{\lambda}} B^*A^*[BA]AB^*[A^*B]\\
                     &=& \lambda^2 B^*A^*[\frac{1}{\lambda} AB]AB^*[\frac{1}{\lambda}BA^*]\\
                     &=&  B^*A^*ABA[B^*B]A^*\\
                       &=&  B^*A^*ABA[BB^*]A^*\\
                           &=&  (AB)^*AB AB (AB)^*
\end{eqnarray*}
We obtain that $AB$ is binormal
\end{pr}
\begin{thm}
Let $A$ be $k$-hyponormal and $B$ isometry such that $AB=\lambda BA\neq 0,\;\; \lambda\in\mathbb{C}$. The following statements are equivalent
\begin{enumerate}
  \item $AB$ is $k$-hyponormal
  \item $\sigma(AB)\neq \{0\}$
  \item $|\lambda|=1$
\end{enumerate}
\end{thm}
\begin{pr}
It is clear $(1)\Rightarrow (2)\Rightarrow (3)$. So we show $(3)\Rightarrow (1)$.
For any unit vector $x\in H$ scince $A$ is $k$-hyponormal and $B$ is isometry then

\begin{eqnarray*}
\|(AB)x\|^k &=& \|A(Bx)\|^k  \\
     &\leq& \|A^k(Bx)\|\|Bx\|^{k-1}  \\
     &\leq& \|A^k(Bx)\|
\end{eqnarray*}

On the other hand

By induction we show that $(AB)^k=\lambda^{\frac{k(k-1)}{2}} B^{k-1}A^kB$ for every $k\in\mathbb{N}^*$.\\
For $k=1$, $(AB)^1= \lambda^{\frac{1(1-1)}{2}} B^{1-1}A^1B$.
Assume that it holds for $k\geq2$, we obtain
\begin{eqnarray*}
(AB)^{k+1}=AB(AB)^k&=&(\lambda BA)(\lambda^{\frac{k(k-1)}{2}} B^{k-1}A^kB)\\
                  &=& \lambda^{\frac{k(k-1)}{2}+1}BA B^{k-1}A^kB\\
                  &=& \lambda^{\frac{k(k-1)}{2}+1}B(AB)B^{k-2}A^kB\\
                    &=& \lambda^{\frac{k(k-1)}{2}+1}B(\lambda BA)B^{k-2}A^kB\\
                     &=& \lambda^{\frac{k(k-1)}{2}+2}B^2 AB^{k-2}A^kB\\
                     &:&\\
                     &:&\\
                     &:&\\
                     &=& \lambda^{\frac{k(k-1)}{2}+k}B^k AB^{k-k}A^kB\\
                     &=& \lambda^{\frac{(k+1)k}{2}}B^k A^{k+1}B
                      \end{eqnarray*}
We conclude that $(AB)^k=\lambda^{\frac{k(k-1)}{2}} B^{k-1}A^kB$ for every $k\in\mathbb{N}^*$.\\
Because $B$ is isometry and $|\lambda|=1$ it follows that
\begin{eqnarray*}
 \|(AB)^k\|  &=& \|\lambda^{\frac{k(k-1)}{2}} B^{k-1}A^kBx\|  \\
       &=& |\lambda|^{\frac{k(k-1)}{2}}\|B^{k-1}A^kBx\|  \\
       &=& \|A^kBx\|
\end{eqnarray*}
Therefore it holds  $\|(AB)x\|^k \leq \|(AB)^kx\|$  for any unit vector $x$ and $AB$ is \\ $k$-hyponormal, which completes the proof.
\end{pr}

\begin{thm}
Let $A,B\in \mathcal{B}(\mathcal{H})$ such that $AB=\lambda BA\neq 0, \;\;\lambda\in\mathbb{C}$.  Then
\begin{enumerate}
  \item if $A^*$ is $M_1$-hyponormal and $B$ is $M_2$-hyponormal, then $|\lambda|\leq (M_1M_2)^\frac{1}{2}$

  \item if $A$ is $M_1$-hyponormal and $B^*$ is $M_2$-hyponormal, then $|\lambda|\geq (M_1M_2)^{-\frac{1}{2}}$

\end{enumerate}

\end{thm}

\begin{pr}
\begin{enumerate}
  \item We have
\begin{eqnarray*}
  |\lambda|\|BA\|&=&\|\lambda BA\|\\
         &=&  \|AB\|\\
       &=& \|B^*A^*AB\|^{\frac{1}{2}}  \>\>\>\>(\|T\|=\|TT^*\|^{\frac{1}{2}}\\
       &\leq & M_1^{\frac{1}{2}}\|B^*AA^*B\|^{\frac{1}{2}} \>\> (A^* \>\textit{is}\> M_1-hyponormal:\>  A^*A\leq M_1AA^*)\\
       &\leq&   M_1^{\frac{1}{2}} \|A^*B\|\>\>(\|T^*T\|^{\frac{1}{2}}=\|T\|)\\
       &\leq& M_1^{\frac{1}{2}}  \|A^*BB^*A\|^{\frac{1}{2}}\>\>\>\>(\|T\|=\|TT^*\|^{\frac{1}{2}})  \\
       &\leq & (M_1M_2)^{\frac{1}{2}} \|A^*B^*BA\|^{\frac{1}{2}} \>\> (B \>\textit{is}\> M_2-hyponormal:\>BB^*\leq M_2B^*B)  \\
       &\leq& (M_1M_2)^{\frac{1}{2}} \|BA\|\>\>\>\>(\|T^*T\|^{\frac{1}{2}}=\|T\|).
\end{eqnarray*}
         Hence $|\lambda|\|BA\| \leq (M_1M_2)^{\frac{1}{2}} \|BA\|$ and $|\lambda|\leq (M_1M_2)^{\frac{1}{2}}$\\
          \item Since $AB=\lambda BA$ and $\lambda\neq 0$ we have $BA=\lambda^{-1}AB$ and apply (1) we obtain $|\lambda^{-1}|\leq (M_2M_1)^{\frac{1}{2}}$,
          therefore  $|\lambda|\geq (M_2M_1)^{-\frac{1}{2}}$.\\
\end{enumerate}
\end{pr}

 \begin{coro}
 Let $A,B\in \mathcal{B}(\mathcal{H})$ such that $AB=\lambda BA\neq 0, \;\;\lambda\in\mathbb{C}$.  Then
\begin{enumerate}
  \item if $A^*$ and $B$ are hyponormal, then $|\lambda|\leq 1$
  \item if $A$ and $B^*$ are hyponormal, then $|\lambda|\geq 1$\\
\end{enumerate}
\end{coro}

\begin{thm}
Let $A,B\in \mathcal{B}(\mathcal{H})$ such that $AB=\lambda BA\neq 0,\;\; \lambda\in\mathbb{C}$.  Then
if $A^*$ is $M_1$-hyponormal and $B$ is $M_2$-hyponormal then $A^*B$ and $BA^*$ are $M_1M_2|\lambda|^2$-hyponormal

\end{thm}

\begin{pr}
We have
\begin{eqnarray*}
   (A^*B)^*A^*B &=& B^*AA^*B  \\
       &\geq &  M_1B^*A^*AB \\
       &\geq &  M_1 \bar{\lambda}A^*B^*\lambda BA \\
         &\geq &  M_1 |\lambda|^2 A^*B^*BA\\
          &\geq &  M_1|\lambda|^2A^*M_2BB^* A\\
           &\geq&  M_1M_2|\lambda|^2(B^* A)^*B^* A
           \end{eqnarray*}

Therefore $A^* B$ is $M_1M_2|\lambda|^2$-hyponormal.\\

Same way we prove that $BA^* $ is $M_1M_2|\lambda|^2$-hyponormal.\\

\end{pr}

\end{document}